\documentclass[12pt]{report}
\usepackage{amsmath}
\usepackage{amsfonts}
\usepackage{amssymb}
\usepackage{mathrsfs}
\usepackage{graphicx}
\usepackage{cancel}
\usepackage[export]{adjustbox}
\usepackage{adjustbox}
\usepackage{lipsum}% http://ctan.org/pkg/lipsum
\usepackage{hyperref}

\DeclareMathAlphabet{\mathpzc}{OT1}{pzc}{m}{it}

\newtheorem{theorem}{Theorem}[chapter]
\newtheorem{corollary}{Corollary}[chapter]

\newtheorem{lemma}{Lemma}[chapter]

\begin{document}

\begin{flushleft}
JOURNAL OF COMBINATORIAL THEORY, Series B \textbf{42}, 146-155 (1987)
\end{flushleft}
$$$$$$$$
\begin{center}
\textbf{\large Circuit Preserving Edge Maps II}

\end{center}
\begin{center}
\textbf{Jon Henry Sanders }
\end{center}
$$$$
\begin{center}
JHS Consulting          jon\_sanders@partech.com\\

\end{center}

\noindent\par In Chapter 1 of this article we prove the following. Let $f:G\rightarrow G^\prime$ be a \textit{circuit surjection,} i.e., a mapping of the edge set of $G$ onto the edge set of $G^\prime$ which maps circuits of $G$ onto circuits of $G^\prime,$ where $G,G^\prime$ are graphs without loops or multiple edges and $G^\prime$ has no isolated vertices. We show that if $G$ is assumed finite and 3-connected, then $f$ is induced by  a vertex isomorphism. If $G$ is assumed 3-connected but not necessarily finite and $G^\prime$ is assumed to not be a circuit, then $f$ is induced by a vertex isomorphism. Examples of circuit surjections $f:G\rightarrow G^\prime$ where $G^\prime$ is a circuit and $G$ is an infinite graph of arbitrarily large connectivity are given. In general if we assume $G$ two-connected and $G^\prime$ not a circuit then any circuit surjection $f:G	\rightarrow G^\prime$ may be written as the composite of three maps, $f(G)=q(h(k(G))),$ where $k$ is a $1-1$ onto edge map which preserves circuits in both directions (the``2-isomorphism'' of Whitney(\textit{Amer. J. Math.} 55(1993), 245-254 ) when $G$ is finite), $h$ is an onto edge \textit{circuit injection} (a 1-1 circuit surjection). Let $f: G\rightarrow M$ be a 1-1 onto mapping of the edges of $G$ onto the cells of $M$ which takes circuits of $G$ onto circuits of $M$ where $G$ is a graph with no isolated vertices, $M$ a matroid. If there exists a circuit $C$ of $M$ which is not the image of a circuit in $G$, we call $f$ \textit{nontrivial}, otherwise \textit{trivial}. In Chapter 2 we show the following. Let $G$ be a graph of even order. Then the   statement        `` no nontrivial map $f: g\rightarrow M$ exists, where $M$ is a binary matroid,'' is equivalent to ``$G$ is Hamiltonian.'' If $G$ is a graph of odd order, then the statement ``no nontrivial map $f:G\rightarrow M$ exists, where $M$ is a binary matroid''  is equivalent to ``$G$ is almost Hamiltonian'', where we define a  graph $G$ of order $n$ to be \textit{almost Hamiltonian} if every subset of vertices of order $n-1$ is contained in some circuit of $G$.

\begin{center}
\textbf{INTRODUCTION AND PRELIMINARY DEFINITIONS}
\end{center}
\noindent\par The results obtained in this paper grew from an attempt to generalize the main theorem of [1]. There it was shown that any \textit{circuit injection} (a 1-1 onto edge map $f$ such that if $C$ is a circuit then $f(C)$ is a circuit from a 3-connected (not necessarily finite)) graph $G$ onto a graph $G^\prime$ is induced by a vertex isomorphism, where $G^\prime$ is assumed to not have any isolated vertices. In the present article we examine the situation when the 1-1 condition is dropped (Chapter 1). An interesting result then is that the theorem remains true for finite (3-connected ) graphs $G$ but not for infinite $G$.
\par In Chapter 2 we retain the 1-1 condition but allow the image of $f$ to be first an arbitrary matroid and second a binary matroid.
\par Throughout this paper we will assume that graphs are undirected without loops or multiple edges and not necessarily finite unless otherwise stated. We will denote the set of edges of a graph $G$ by $E(G)$ and  the set of vertices of $G$ by $V(G)$. We will also use the notation $G=(V,E)$ to indicate $V=V(G), E=E(G)$ when $G$ is a graph. The graph $G: A$ will be the graph with edge set $A$ and vertex set $V(G)$. The abuse of language of referring to a set of edges $S$ as a graph (usually a subgraph of a given graph) will be tolerated where it is understood that the set of vertices of such a graph is simply the set of all vertices adjacent to any edge of $S$.
\par A subgraph $P$ of a graph $G$ is a \textit{suspended chain} of $G$ if $|V|\geq 3, |V|$ finite and there exists two distinct vertices $v_1,v_2\in V$, the endpoints of $P$ such that $\deg_Pv_1=1,~~ \deg_P v_2=1$, and $\deg_Pv=\deg_Gv=2$ for $v\in V, v\neq v_1,v_2$, where $V=V(P).$ We shall also refer to the set of edges of $P$ as a suspended chain. The notation $\mathscr{C}(v)$ will be used to indicated the set of edges adjacent to the vertex $v$ in a given graph.
\par A \textit{circuit surjection f} of $G$ onto $G^\prime,$ denoted by $f:G\rightarrow G^\prime$, is an onto map of the edge set of $G$ onto the edge set of $G^\prime$ such that if $C$ is a circuit of $G$ then $f(C)$ is a circuit of $G^\prime$. We also understand the terminology $f:G\rightarrow G^\prime$ is a circuit surjrction to preclude the possibility of $G^\prime$ having isolated vertices.
\chapter{\small 1. CIRCUIT SURJECTIONS ONTO GRAPHS}
\numberwithin{theorem}{chapter}

\begin{lemma}
Let $f:G\rightarrow G^\prime$ be circuit surjection where $G$ is 2-connected and $G^\prime$ is not a circuit. Let $e$ be an edge of $G^\prime$. Then if $C$ is circuit of $G$ such that $C$ contains at least one element of $f^{-1}(e)$ then $C$ contains every element of $f^{-1}(e).$
\end{lemma}
\noindent\par\textit{Proof.}~~~~ First, we note that $G^\prime$ is 2-connected since if $e_1,e_2$ are two distinct edges of $G^\prime$ then $f(C)$ is a circuit which contains $e_1$ and $e_2$ where $C$ is any circuit of $G$ which contains $h_1,h_2$ such that $h_1\in f^{-1}(e_1),h_2\in f^{-1}(e_2).$ Let $v_1,v_2$ be the vertices adjacent to $e$. Let $P(v_1,v^\prime)$ be a path in $G^\prime$ of minimal length such that $v^\prime$ is a vertex of degree greater than 2. Define $S=\mathscr{C}(v^\prime)-\{h\}$ if $v^\prime\neq v_1,~~ S=\mathscr{C}(v^\prime)-\{e\}$ if $v^\prime=v_1,$ where $h$ is the edge in $P(v_1,v^\prime)$ adjacent to $v^\prime.$\\

FACT 1.~~ Any circuit of $G^\prime$ which contains $e$ must contain one and only one element of $S$.
\par Let $a_\alpha,\alpha\in I$ be the elements of $S$ and let $A=f^{-1}(e),A_\alpha=f^{-1}(a_\alpha),\alpha\in I.$ Then Fact 1 implies\\

FACT 2.~~ If $C\cap A\neq\O$ for $C$ a circuit of $G$ then $C\cap A_\alpha\neq\O$  is true for one and only $\alpha\in I.$
\par Let $C_0$ be a circuit which contains an edge of $A$. We will show that the assumption $C_0\not\supset A$ leads to a contradiction of Fact2. Denote by $B$ the unique set $A_{\alpha 0}, \alpha_0\in I$ such that $C_0\cap A_{\alpha 0}\neq\O.$ Let $D=A_{\alpha 1}, \alpha_1\neq \alpha_0$ (since $|I|=|S|\geq 2$, this is possible) and let $d\in D$. Since $G$ is 2-connected and $d\notin C_0$ there is a path $P_3(q_0,q_1),d\in P_3(q_0,q_1)$ where $q_0,q_1$ are distinct vertices of $C_0$ and $P_3(q_0,q_1)$ is edge disjoint from $C_0.$ Denote by $P_1(q_0,q_1)$ and $P_2(q_0,q_1)$ the two paths such that $C_0=P_1(q_0,q_1)\cup P_2(q_0,q_1).$ Now $P_i\cap A\neq\O$ and $P_i\cap B\neq\O$ is not possible, $i=1$ or $2$, since then $P_3\cup P_i$ would be a circuit which violates Fact 2. Thus $P_i\cap A\neq\O (P_i\cap B=\O),$ and $P_j\cap B\neq\O (P_j\cap A=\O)$ where either $i=1, j=2,$ or $j=1, i=2,$ say, the former (Fig. 1).
\par Suppose now there exists an edge $k\in A,k\notin C_0$. Now $k\in P_3$ is impossible since if that were the case then $P_3\cup P_2$ would be a circuit which violates Fact 2. Thus $k$ is edge disjoint from $G^{\prime\prime}$, where $G^{\prime\prime}$ is the subgraph of $G$ consisting of $P_3\cup P_1\cup P_2.$ Since $G$ is 2-connected there exists a path $P_4(t_0,t_1)$ in $G$ such that $k\in P_4(t_0,t_1),t_0,t_1$ are distinct vertices of $G^{\prime\prime}$ and $P_4(t_0,t_1)$ is edge disjoint from $G^{\prime\prime}$. We now show that no matter where $t_0,t_1$ fall on $G^{\prime\prime}$ a contradiction to Fact 2 arises. For if $G^{\prime\prime}$ has a $t_0-t_1$ path $P_5$ disjoint from $B\cup D,$ then $P_4\cup P_5$ is a circuit intersecting $A$ and hence $P_4$ intersects some $A_\alpha$. Since $P_4$ can be extended to a circuit intersecting $B$ (resp. $D$) this contradicts Fact 2. If $G^{\prime\prime}$ has no such path $P_5$, then it has a $t_0-t_1$ path intersecting both $B$ and $D$ and that path union $P_4$ contradicts Fact 2.\\

\begin{theorem}
Let $f:G\rightarrow H$ be a circuit surjection, where $G$ is 2-connected and $H$ is not a circuit. Then $f$ is the composite of three maps $f(G)=g(h(k(G))),$ where $k$ is a 1-1 onto edge map which preserves circuits in both directions (a ``2-isomorphism'' of [8] when $G$ is finite), $h$ is an onto edge map obtained by replacing suspended chains by single edges (which preserves circuits in both directions) and $q$ is a circuit injection. 
\end{theorem}
$$$$$$$$$$$$$$$$$$$$
%\begin{figure}[hp]
	%\centering
	%\begin{center}
		\includegraphics[scale=1.0]{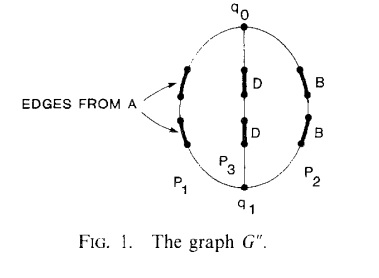}
	%\end{center}
	%\label{fig:fig1}
%\end{figure}

\noindent\par We note that the theorem implies that $f^{-1}(e)$ is a finite set for each edge of $H$  and thus $H$ must be infinite if $G$ is infinite.
\par Theorem 1.1 follows from the fact that (by Lemma 1.1) for any $e\in H$, any two edges of $f^{-1}(e)$ form a minimal cut set (cocycle) It is apparent that $f^{-1}(e)$ can thus be transformed into a suspended chain by a sequence of 2-switchings. This establishes Theorem 1.1 for finite $G$. Theorem 1.1 also holds for infinite $G$ by the same method used in Theorem 4.1 of [3] (where Whitney's 2-isomorphism theorem [8] is extended to the infinite case).
\begin{theorem}
Let $f:G\rightarrow G^\prime$ be a circuit surjection, where $G$ is finite and 3-connected. Then $f$ is induced by a vertex isomorphism.
\end{theorem}
\noindent\par\textit{Proof.}~~~~ We will show that $G^\prime$ cannot be a circuit. For assume $G^\prime$ is a $k$-circuit, $k\geq3.$ Write $G=(V,E)$ and $|V|=n$. Now $f^{-1}(G^\prime-\{e\})$ contains no circuit and thus $|f^{-1}(G^\prime-\{e_i\})|<n,i=1,\ldots,k,$ ~~ where $e_1,\ldots,e_k$ are the edges of $G^\prime.$ But each of $G,$ i.e., each element of $E$ occurs in exactly $k-1$ of the $k$ sets $f^{-1}(G^\prime-\{e_i\},i=1),\ldots,k,$ and $E=\displaystyle \bigcup_{i=1,\ldots,k} f^{-1}(G^\prime-\{e_i\}).$ Thus $(k-1)|E|<kn,$ or $|E|<(k/(k-1))n,$ and thus $|E|<\frac{3}{2}n.$ But $|E|\geq\frac{3}{2}n$ for any (finite) graph each vertex of which is of degree three or greater and thus for any 3-connected finite graph,$\Rightarrow\Leftarrow$. Thus $G^\prime$ cannot be a circuit. Theorem 1.1 thus implies that $f$ is 1-1 so the result follows from [1].
\begin{theorem}
Let $f:G\rightarrow G^\prime$ be a circuit surjection, where $G$ is 3-connected, not necessarily finite and $G^\prime$ is not a circuit. Then $f$ is induced by a vertex isomorphism.
\end{theorem}
\noindent\par\textit{Proof.}~~~~ Theorem 1.1 implies that $f$ must be a 1-1 map so the result follows from [1].\\
\textit{Construction}\\
\noindent\par An $n$-connected graph which has a circuit surjection onto a 3-circuit may be obtained from a sequence of disjoint 2-way infinite paths $P_1,P_2,\ldots,$ such that each vertex of $P_i$ is ``connected'' to $P_{i+1}$ by a tree as indicated in Fig. 2 for $n=4$. (The mapping which takes each edge labeled $i$ onto $e_i,i=1,2,3,$ defines the circuit surjection onto the 3-circuit with edges $e_1,e_2,$ and $e_3$)
\chapter{\small 2. CIRCUIT INJECTIONS ONTO MATROIDS}
\textit{Terminology and Notation}\\
\noindent\par A \textit{matroid} $M$ is an ordered pair of sets $\{S,\mathscr{C}\},$ where $S\neq\O,\mathscr{C}\subseteq 2^S$, which satisfies the following two axioms. Axiom I. $A,B\in\mathscr{C}, A\subseteq B$ implies $A=B$. Axiom II. $A,B\in\mathscr{C}, a\in A\cap B,~~ b\in(A\cup B)-(A\cap B)$ implies there
\newpage
%$$$$$$$$$$$$$$$$$$$$$$$$$$$$$$$$$$$$$$$$$$$$$$$$$$$$$$$$$$$$
%\begin{center}
	\includegraphics{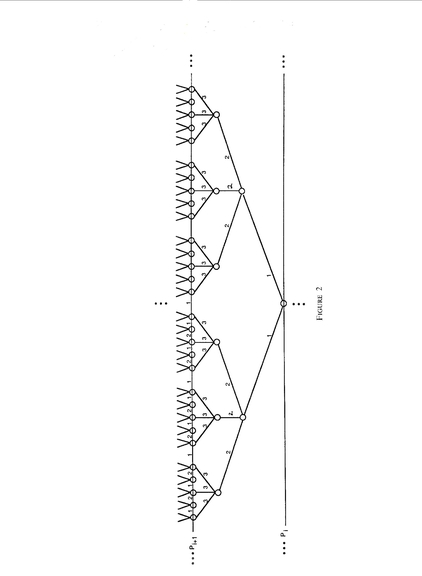}
%\end{center}%Image is here  

%There is an IMAGE HERE
%\newpage
 %%\begin{figure}[hp]
	 %%\centering
		 %\includegraphics[max size={\textwidth}{\textheight}]{Images/fig2.jpg}
	 %%\label{fig:Figure22}
 %%\end{figure}
%\MyincludeGraphics{Images/fig2.jpg}
exists $D\in\mathscr{C}$ such that $D\subseteq A\cup B,a\notin D, b\in D.$ The elements of $S$ are called the cells of $M$, the elements of $\mathscr{C}$ are called the circuits of $M.$
\par The matroid associated with a graph $G_M$, is the matroid whose cells are the edges of $G$ and whose circuits are the circuits of $G$.
\par Let $M=\{S,\mathscr{C}\}, M^\prime=\{S^\prime,\mathscr{C}^\prime\}$ be matroids, and let $f:S\rightarrow S^\prime$ be a 1-1 onto map such that $f(A)\in\mathscr{C}^\prime$ whenever $A\in\mathscr{C}$. Such an $f$ is called a circuit injection of $M$ onto $M^\prime$ denoted by $f:M\rightarrow M^\prime.$ The circuit injection injection $f$ is called nontrivial if there exists $B\in\mathscr{C}^\prime$ such that $B\neq f(A)$ for all $A\in\mathscr{C}.$
\par We can assume without loss of generality that $S=S^\prime, f$ is the identity map and $\mathscr{C}\subseteq\mathscr{C}^\prime$ for a circuit injection $f$. Then $f$ is nontrivial if $\mathscr{C}$ is properly contained in $\mathscr{C}^\prime$.
\par We denote by $A\oplus B$ the $mod$ 2 $addition$ of set $A$ and $B$ which is defined to be the set $(A\cup B)-(a\cap B).$
\par A matroid $(S,\mathscr{C})$ is a binary matroid if for all $A,B\in\mathscr{C}, A\oplus B =\displaystyle \bigcup_{i=1}^k C_i$ for $C_i\in \mathscr{C}, i=1,\ldots,k,~~ C_i\cap C_j=\O, i\neq j,1\leq i, j\leq k$. Given a set $S$ and an arbitrary set $\mathscr{C}\subseteq 2^S$ we denote by $<\mathscr{C}>$ the collection of all sets $A$ such that there exists $k\geq1, C_1,\ldots, C_k\in\mathscr{C}$ and $A=C_1\oplus\cdots\oplus C_k.$
\par We denote by $<\mathscr{C}>_{\min}$ the minimal elements of $<\mathscr{C}>,$~~i.e., the elements $A\in<\mathscr{C}>$ such that $B\in<\mathscr{C}>, B\subseteq A\Rightarrow B=A$. A useful theorem of matroid theory [5, Sects. 1 and 5.3] is that $\{S,<\mathscr{C}>_{\min}\}$ is a binary matroid for arbitrary $\mathscr{C}\subseteq 2^S.$ 
\par We denote the rank of a matroid by $r(M)$. If $A\in\mathscr{C}$ exists such that $|A|=r(M)+1$ we call $A$ a $Hamiltonian$  $circuit$ of $M$, and we call $M$ $Hamiltonian$.\\

\textit{Condition for Trivial/Nontrivial Circuit Injections}\\
\par We would like to establish conditions on a graph $G$ such that all circuit injections $f:G_M\rightarrow N$ are trivial, where $N$ is first assumed to be an arbitrary matroid and second assumed to be a binary matroid. (We note that if $N$ is assumed to be a graphic matroid, i.e., $N=G^\prime_M$ for some graph $G^\prime$ then the theorem of [1] implies that $G$ 3-connected is a condition when ensures no nontrivial circuit injection exists).
\par Since the addition of an isolated vertex to a graph $G$ has no effect on $G_M$ we assume (without loss of generality) that $G$ has no isolated vertices throughout this section to simplify the statements of the theorems.\\
\textit{Remark.}~~~~ The fact that if $M$ is a Hamiltonian matroid (or in particular $G_M$, where $G$ is a Hamiltonian graph) then the only circuit injections $f: M\rightarrow M^\prime$ are trivial, where $M^\prime$ is an arbitrary matroid follows from the fact that $r(M^\prime)=r(M)$ in this case. The converse is also easily established as follows. 

\begin{theorem}
If $G$ is a non-Hamiltonian matroid (or in particular the matroid associated with a graph without Hamiltonian circuits) there exists a nontrivial circuit injection $f:G\rightarrow M$, where $M$ is a (not in general binary) matroid.
\end{theorem}
\noindent\par\textit{Proof.}~~~~ Let the cells of $M$ be the edges of $G$; let the circuits of $M$ be $\mathscr{C}\cup\mathscr{L},$ where $\mathscr{C}$ is the set of circuits of $G$ and $\mathscr{L}$ is the set of all bases of $G$, and let $f$ be the identity map. Then $f$ is a nontrivial circuit injection (the matroid $M$ is the so-called truncation of $G$ see [7]).
\par \textit{Remark}.~~~~ Since matroids of arbitrarily large connectivity exist without Hamiltonian circuits (the duals of complete graphs are one example \footnote{We take the definition of connectivity for matroids from [4, 6]. A property of this definition is that the connectivity of a matroid equals the connectivity of its dual and also the connectivity of the matroid $G^n_M$ associated with the complete graph on n vertices $G^n$ approaches $\infty $ as $n \rightarrow \infty $. Thus the duals of the complete graphs have arbitrarily large connectivity.}) there is no general matroid analogue to the result of [1]. We note that $M$  is never a binary matroid in the construction of Theorem 2.1.
\par A more interesting result is obtained when we restrict $M$ to be an arbitrary matroid, $G$ a graphic matroid.\\

\par DEFINITION.~~~~Let the order of a graph $G$ be $n$. We say $G$ is almost Hamiltonian if every subset of $n-1$ vertices is contained in a circuit.

\begin{theorem}
Let the order of $G$ be even. Then ``no nontrivial circuit injection $f$ exists, $f:G_M\rightarrow B,$ where $B$ is binary'' is true iff $G$ is Hamiltonian. Let the order of $G$ be odd. Then ``no nontrivial circuit injection $f:G_M\rightarrow B$ exists where $B$ is binary'' is true iff $G$ is almost Hamiltonian.
\end{theorem}
We abbreviate ``no nontrivial circuit injection $f:G_M\rightarrow B$ exists, where $B$ is binary'' by saying ``$G$ has no nontrivial map.'' To prove the theorem we need the following
\begin{lemma}
$G$ has no nontrivial map implies ``if $v_1,\ldots,v_n$ are vertices of odd degree in $S$, for any subgraph $S$ of $G$, then there exists a circuit $C$ of $G$ such that $v_1,\ldots,v_n$ are vertices of C.''
\end{lemma}
\noindent\par\textit{Proof.}~~~~ Let $\mathscr{C}$ be the set of circuits of $G$, $S$ a subset of edges of $G$. Let $\mathscr{C}^\prime=<\mathscr{C}\cup\{S\}>_{\min}$. Then $f:\{E,\mathscr{C}\}\rightarrow\{E,\mathscr{C}^\prime\},$  where $f$ is the identity map, will be a circuit injection unless $\mathscr{C}\not\subseteq\mathscr{C}^\prime$, i.e., unless there exists $A\in<\mathscr{C}\cup\{S\}>_{\min},C\in\mathscr{C}$ and $A$ is properly contained in $C$, i.e., unless
\begin{equation}
S\oplus C_1\oplus\cdots\oplus C_k\subset C	\qquad\qquad \mbox{   for   }\qquad C_i\in\mathscr{C}, i=1,\ldots,k.
\end{equation}
Now if $S$ has a vertex $v$ of odd degree in $S$ then $\mathscr{C}\neq<\mathscr{C}\cup\{S\}>_{\min}$ so $f$ will be a nontrivial circuit injection unless $(2.1)$ holds. But $v$ of odd degree in $S$ implies $v$ will be of odd degree in $S\oplus C_1\oplus\cdots\oplus C_k$ and thus $v$ must be contained in $C$. (If vertex $q$ is of even degree in $S$ then all edges adjacent to it could cancel in $S\oplus C_1\oplus\cdots\oplus C_k$ and thus $q\notin C$ is possible).
\begin{corollary}
$G$ has no nontrivial map implies $G$ is 2-connected.
\end{corollary}
\noindent\par\textit{Proof.}~~~~ We show given $q_1\neq q_2,$ vertices of $G$, there exist $C\in\mathscr{C}$ with $q_1,q_2$ vertices of $C$. First assume there exists the edge $e=(q_1,q_2)$ in $G$. Then taking $S=\{e\}$ in the hypothesis of Lemma 2.3 yields $C$. Otherwise choose an edge a adjacent to $q_1$ and an edge $b$ adjacent to $q_2$ (since $G$ has no isolated vertices this is possible) and put $S=\{a,b\}$ to get $C$.
\par We prove the implications of Theorem 2.2 separately in the following two lemmas.
\begin{lemma}
$|G|=2N$ and $G$ has no nontrivial map $\Rightarrow G$ is Hamiltonian; $|G|=2N+1$ and $G$ has no nontrivial map $\Rightarrow G$ is almost Hamiltonian.
\end{lemma} 
\noindent\par\textit{Proof.}~~~~ Let $C$ be a circuit of $G$ and let $G$ have no nontrivial map, $|G|$ odd or even.
\par FACT 1.~~ If $C$ is even and there exist two distinct vertices $v_1,v_2$ of $G$ not on $C$ then $C$ is not of maximal order.
\par \textit{Proof of Fact 1.}~~ Let $q_1,q_2$ be two distinct vertices of $C$. Then by Menger's Theorem (since $G$ is 2-connected ) there exists a pair of vertex disjoint paths $P(v_1,q_1), P(v_2,q_2)$ or $P(v_1,q_2), P(v_2,q_1)$. In either case there exists a pair of distinct vertices $v_1^\prime,v_2^\prime$ not on $C$ such that $(v_1^\prime, q_1),(v_2^\prime,q_2)$ are edges of $G$. If $q_1,q_2$ are separated by an odd (even) number of edges  in $C$ there exists a subgraph of $G$ having $|C|+2$ odd vertices as in Fig. 3(A) (3(B)) and thus $C$ is not maximal by Lemma 2.1.
%$$$$$$$$$$$$$$$$$$$$$$$$$$$$$$$$$$$$$$$$$$$$$$$$$$$$$$$$
%\begin{center}
	\newpage
	%\begin{flushleft}
	
\begin{center}
		\includegraphics[width=1\textwidth]{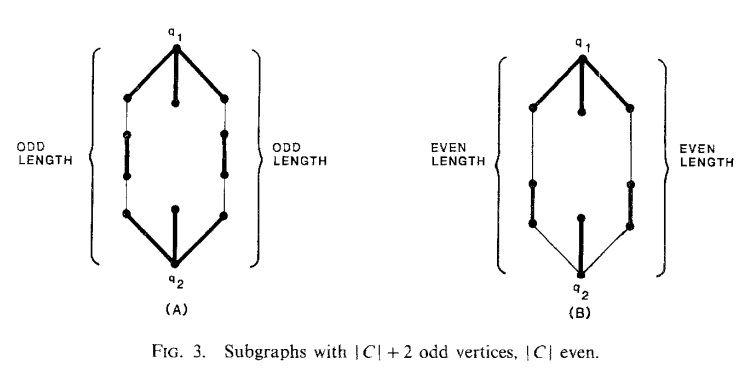}
\end{center}
	%\
	%\end{flushleft}
%end{center}%Image is here  
%Image is here
$$$$$$$$

	\begin{flushright}
		\includegraphics{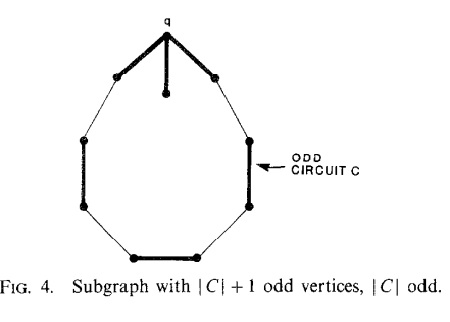}
	\end{flushright}

\par FACT 2.~~ If $|C|$ is odd and there exists a vertex $v_1\in G$ not on $C$ then $C$ is not maximal.
\par \textit{Proof of Fact 2.}~~ By the connectivity of $G$ we have $(v_1,q)$ is an edge for some vertex $q$ on $C$. We construct a subgraph having $|C|+1$ odd vertices as in Fig. 4 and apply Lemma 2.1.
\par If $|G|=2N$, Facts 1 and 2 imply that a circuit of maximal length is a Hamiltonian circuit. If $|G|=2N+1$, Facts 1 and 2 imply either $G$ is Hamiltonian (in which case it is also almost Hamiltonian) or a maximal circuit is of length $2N$. Let $C$ be a circuit of length $2N,v$ the vertex of $G$ not on $C, q$ a vertex on $C$ such that $(v,q)$ is an edge. We can find a subgraph of $G$ all vertices of which are of odd degree containing $v$ and all other vertices of $C$ other than an arbitrary vertex $v^\prime$ of $C$ as in Fig. 5. Thus $G$ is almost Hamiltonian by Lemma 2.1.
\begin{lemma}
Let $G$ be an almost Hamiltonian graph, $|G|=2N+1$. Then $G$ has no nontrivial map. Let $G$ be a Hamiltonian graph, $|G|=2N$. Then $G$ has no nontrivial map.
\end{lemma}
\noindent \par \textit{Proof.}~~~~ Case 1. $|G|=2N+1.$ Suppose otherwise, i.e., let $f:(E,\mathscr{C})\rightarrow(E,\mathscr{C}^\prime)$ be a nontrivial circuit injection, where $E$ are the edges of $G,\mathscr{C}$ are the circuits of G, and $\mathscr{C}^\prime$ properly contains $\mathscr{C}.$ Let $C$ be a circuit of $G$, \\
% There is Image here 
%$$$$$$$$$$$$$$$$$$$$
\newpage
\includegraphics[width=1\textwidth]{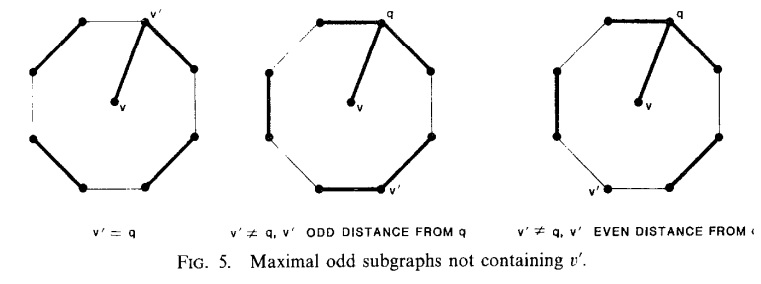}

$$$$

$|C|=2N, q$ a vertex of $G$ not on $C,e^\prime$ an edge of $G$ adjacent to $q$ and some vertex $v$ of $C$, and $e$ an edge of $C$ adjacent to $v$.
\par Then $P=(C-\{e\})\cup\{e^\prime\}$ is a Hamiltonian path of $G$ (i.e., a path which contains every vertex) and $P$ is a dependent set of $\{E,\mathscr{C}^\prime\}$(since otherwise $r(E,\mathscr{C})=r(E,\mathscr{C}^\prime)=2N$ and $f$ must be trivial). Let $T\in\mathscr{C}^\prime, T\notin\mathscr{C}, T\subseteq P.$ Now $T$ has at most $2N$ odd vertices, $v_1,\ldots,v_s$, since the sum of the degrees of all the vertices of $T$ is even and $T$ has at most $2N+1$ vertices. Let $C^\prime$ be  a circuit of $G$  which contains $v_1,\ldots,v_s$. Let $T\subseteq T$ be the set of edges of $T$ not contained in $C^\prime.$ Then $T^\prime\subseteq P$ is the union of vertex disjoint paths $P_1,\ldots,P_k$ and the endpoints $b_i,e_i$ of $P_i$ are on $C^\prime$.Let $C_i^\prime$ be one of the two paths in $C^\prime$ with endpoints $b_i,e_i$ of $P_i$ are on $C^\prime$. Let $C_i^\prime$ be one of the two paths in $C^\prime$ with endpoints $b_i,e_i$ and define $k$ circuits of $G$ by $C_i=C_i^\prime\cup P_i, i=1,\ldots,k.$ Then $T\oplus C_i\oplus\cdots\oplus C_k\subseteq C^\prime$ contradicting the definition of $T$.\\
\par Case 2.~~ $|G|=2N$. If $G$ is Hamiltonian of arbitrary order then $G$ has no nontrivial map as noted in an earlier remark.
\par Lemmas 2.2 and 2.3 establish Theorem 2.2. The existence of almost Hamiltonian graphs of odd order which are not Hamiltonian is shown in [2]. Thus there are graphs which are not Hamiltonian for which no nontrivial map exists.
\par\textit{Remark}.~~ The duals of the matroids of complete graphs of order 5 or more provide a counter example to the assertion that an $n$ exists such that if a binary matroid $M$ has a connectivity $n$ no nontrivial map $f:M\rightarrow M^\prime$ exists, where $M^\prime$ is a binary matroid. For if $G_n$ is the complete graph of $n$ vertices let $M_n^\prime=<B_n\cup\{E_n\}>_{\min},$ where $E_n=E(G_n)$ and $B_n$ is the set of bonds of $G_n.$ Then $f:M_n\rightarrow M_n^\prime,$ where $M_n$ is the dual of $G_n$, and $f$ is the identity map, is a nontrivial map, since $a\oplus E_n\not\subset b$ for $a,b\in B_n$ when $n\geq 5$ and $a_1\oplus\cdots\oplus a_k$  where $a_i\in B_n,$ $1\leq i \leq k.$
\section*{\centering\small ACKNOWLEDGMENT}
The author would like to thank the referee for many helpful suggestions. 
$$$$
\section*{\centering\small REFERENCES}
\begin{description}
\item 1. J.H. SANDERS AND D. SANDERS, Circuit preserving edge maps,\textit{J. Combin. Theory Ser}. B \textbf{22} (1977),91-96.
\item 2. C. THOMASSEN, Planner and infinite hypohamiltonian and hypotraceable graphs, \textit{Discrete math.} \textbf{14 } (1976),377-389.
\item 3.  C. THOMASSEN, Duality of infinite graphs, \textit{J. Combin. Theory Ser.} B \textbf{33} (1982), 137-160.
\item 4. W.T. TUTTE, Menger's theorem for matroids, \textit{J. Res. Nat. Bur. Standards} B \textbf{69}(1964, 49-53).
\item 5. W.T. TUTTE, Lectures on matroids,\textit{J. Res. Nat. Bur. Standards} B \textbf{68}(1965),1-47.
\item 6. W.T. TUTTE, Connectivity in matroids, \textit{Canad. J. Math.} \textbf{18} (1966), 1301-1324
\item 7. D.J. A. WELSH, ``Matroid Theory,'' Academic Press, London/ New York, 1976.
\item 8. H. WHITNEY, 2- isomorphic graphs, \textit{Amer. J. Math.} \textbf{55}(1933), 245-254. 
\end{description}
%$\wp \mathscr{S} \mathscr{C}$

\end{document}